\documentclass[12pt]{amsart}
\usepackage{amsmath,amscd,amsthm,amsfonts, amssymb,amsxtra}
\usepackage[all]{xy}
\SelectTips{cm}{}
\newtheorem{thm}{Theorem}[section]  
\newtheorem*{un-no-thm}{Theorem}
\newtheorem{lem}[thm]{Lemma}         
\newtheorem{prop}[thm]{Proposition}  

\newtheorem{bigthm}{Theorem}

\theoremstyle{definition}
\newtheorem{defn}[thm]{Definition}   

\theoremstyle{definition}
\theoremstyle{definition}

\theoremstyle{remark}
\newtheorem{rem}[thm]{Remark}
  \theoremstyle{remark}
 \newtheorem{rems}[thm]{Remarks}

\newtheorem*{ack}{Acknowledgement}
 
\newtheorem*{out}{Outline}

\begin{document}
\title{The Dualizing Spectrum, II.}
\date{\today}
\author{John R. Klein}
\address{Wayne State University, Detroit, MI 48202}
\email{klein@math.wayne.edu}
\begin{abstract} To an inclusion
$H \subset G$ of topological groups, 
we associate a spectrum $D_{H\subset G}$, which
coincides with the
dualizing spectrum $D_G$ of \cite{Klein_dualizing}
when $H=G$. We also introduce a fibered spectrum analogue.

The main application 
is to give a purely homotopy theoretic construction
of Poincar\'e embeddings in stable codimension.
\end{abstract}
\thanks{The author is partially supported by the NSF}
\maketitle
\setlength{\parindent}{15pt}
\setlength{\parskip}{1pt plus 0pt minus 1pt}
\def\Top{\bold T\bold o \bold p}
\def\wTop{\text{\rm w}\bold T}
\def\wT{\text{\rm w}\bold T}
\def\Sp{\bold S\bold p}
\def\vo{\varOmega}
\def\vs{\varSigma}
\def\smsh{\wedge}
\def\flush{\flushpar}
\def\id{\text{id}}
\def\dbslash{/\!\! /}
\def\codim{\text{\rm codim\,}}
\def\:{\colon}
\def\holim{\text{\rm holim\,}}
\def\hocolim{\text{\rm hocolim\,}}
\def\hodim{\text{\rm hodim\,}}
\def\hocodim{\text{hocodim\,}}
\def\Bbb{\mathbb}
\def\bold{\mathbf}
\def\Aut{\text{\rm Aut}}
\def\cal{\mathcal}
\def\Sec{\text{\rm sec}}
\def\Secst{\text{\rm sec}^{\text{\rm st}}}
\def\maps{\text{\rm map}}
\setcounter{tocdepth}{1}

\section{Introduction \label{intro}}
An embedding of manifolds $P \subset N$ and a compact tubular neighborhood
$D$ of $P$ in $N$ gives rise to a stratification 
$$ 
N = D \cup_{\partial D} C \, ,
$$
where $C$ is the closure of $N-D$.
From a homotopy theorist's perspective, it is legitimate
to consider a variant of the above in which the manifolds
are replaced by Poincar\'e duality spaces, and equality 
is replaced by homotopy equivalence. This yields 
the notion of {\it Poincar\'e embedding.} It turns out that the
`normal datum' $\partial D \to D$ has the homotopy type of
a spherical fibration over $P$.

One the early triumphs 
of surgery theory is the Browder-Casson-Sullivan-Wall
theorem.  It says that a Poincar\'e embedding of PL manifolds in 
codimension at least three can always be lifted 
to a piecewise linear embedding. Thus the problem of 
finding manifold embeddings in these codimensions 
is reduced to a problem in homotopy theory. Even so,
at least one golden era maven opined that
the resulting homotopy theoretic problems
seem ``in some respects to be harder than the original geometrical 
problems'' \cite[p.\ 119]{Wall-book}. 
 
Fortunately, much of the homotopy theory was tractable,
although it took a quarter century to realize it.
By the late 1990s, 
an array of Poincar\'e embedding
results were proved by homotopy theoretic means 
(see e.g., \cite{Richter}, \cite{haef}, \cite{haef2},
\cite{Klein_compression}, \cite{Klein_dualizing}, \cite{Lambrechts-Stanley},
\cite{LSV}; some of these works address a program begun by 
B.\ Williams in the late 1970s \cite{Wi1},\cite{Wi2}).

Actually, the first Poincar\'e embedding result in the literature predates the 
Browder-Casson-Sullivan-Wall theorem: 
Spivak's thesis, somewhat reformulated, 
delivered Poincar\'e embeddings in euclidean space 
in large codimension  and showed that such embeddings
are unique up to concordance \cite{Spivak}. 
The normal datum in this case is called the {\it Spivak fibration}. 

Spivak's results were proved manifold theoretically.
A few years ago, Bill Dwyer and the author (independently)
observed how the Spivak fibration could be obtained by an entirely
homotopy theoretic procedure \cite{Klein_dualizing}.\footnote{The construction
has recently been used to prove that finite  loop spaces have the homotopy type
of smooth manifolds, 
partially settling a question of Browder \cite{finite-loop}.}

The main purpose of this paper is to extend 
some ideas of \cite{Klein_dualizing}, with the aim of giving
a homotopy theoretic proof of the following theorem:

\begin{bigthm}[Existence] \label{embed_thm} Let $f\: P \to N$ be a map
of finitely dominated Poincar\'e duality spaces, where $P$ is without boundary.
Then there is an integer $j > 0$ and
a Poincar\'e embedding of the composite map
$$
\begin{CD} 
f_j\:P @> f >>N @>\subset >> N\times D^j 
\end{CD}
$$
\end{bigthm}

\begin{rem} In \cite{haef}, we used this result as the basis
step of an downward induction on codimension to produce
Poincar\'e embeddings in the range where, roughly, ``connectivity exceeds 
generic double point dimension.'' The proof we gave in \S3 of that paper
 used a manifold thickening of the ambient space together with transversality.

In \cite{haef} we considered the more general situation in which $P$ is merely
a finite complex.
The methods of the current paper do extend to this more general case, 
but for reasons of exposition we decided to retain the
assumption that $P$ satisfies Poincar\'e duality. 
\end{rem}

We shall also prove the following uniqueness result:

\begin{bigthm}[Uniqueness]\label{uniqueness} If $j$ is sufficiently large, then
any two Poincar\'e embeddings of $f_j\:P\to N\times D^j$ are concordant.
\end{bigthm}

\begin{out} \S2 is language. In \S3 we recall the dualizing spectrum
of a topological group and extend it to inclusions. In \S4 we prove
Theorem A in the connected case. In \S5 we describe a fibered spectrum
which has the homotopy type of the unreduced Borel construction of
the dualizing spectrum. In \S6 we complete the proof of Theorem A.
In \S7 we prove Theorem B. In \S8 we present a type of fiberwise
duality, due to Bill Richter, which used in the proof of Theorem B.
\end{out}

\begin{ack} I am indebted to Bill Richter for explaining to me
the duality theory described in \S8.
\end{ack}

\section{Language}

This section  somewhat abridged. The reader
is directed to \cite[\S2]{Klein_dualizing} for a more complete
treatment.

All spaces below will be compactly generated.
A {\it weak equivalence} of spaces means a weak homotopy equivalence.
A space is {\it homotopy finite} it is weak equivalent to
a finite cell complex. It is {\it finitely dominated} if
it is a homotopically a retract of a finite cell complex.

\subsection*{Poincar\'e spaces} 
A finitely dominated space $X$ is a {\it Poincar\'e duality space} 
of (formal) dimension $d$ if there exists 
a bundle of coefficients $\mathcal L$ which is locally isomorphic
to $\Bbb Z$, and a fundamental class $[X] \in H_d(X;{\mathcal L})$
such that the associated cap product homomorphism
$$
\cap [X]\:H^*(X;M) \to H_{d{-}*}(X;{\mathcal L}\otimes M)
$$
is an isomorphism in all degrees. Here, $M$ denotes any 
bundle of coefficients (cf.\  \cite{Wall_PD}, \cite{Klein_PD}). 

Similarly, one has the definition of Poincar\'e space
$X$ {\it with boundary} ${\partial X}$ (also called a
{\it Poincar\'e pair} $(X,\partial X)$).
Here, one assumes both $X$ and ${\partial X}$ are finitely
dominated and there is a fundamental class 
$[X] \in H_d(X,\partial X;{\mathcal L})$
such that 
$$
\cap [X]\:H^*(X;M) \to H_{d{-}*}(X,\partial X; {\mathcal L}\otimes M)
$$
is an isomorphism.  Additionally, if $[\partial X]$ is
 the image of $[X]$ under the boundary homomorphism
$H_d(X,\partial X;{\cal L}) \to 
H_{d-1}(\partial X;{\cal L}_{|\partial X})$, one 
also requires 
$$
({\cal L}_{|\partial X},[\partial X])
$$ 
to equip $\partial X$ with the structure of a Poincar\'e space. 

The above definition makes sense even if the map ${\partial X} \to X$ fails
to be an inclusion, since one can always take a mapping cylinder to convert the
map into a cofibration. We shall sometimes work in this more general
context.

\subsection*{Poincar\'e embeddings}
Let $P^p$ and $N^n$ be a Poincar\'e spaces 
of respective dimensions $p$ and $n$,
with ${\partial P} = \emptyset$ 
(but where $N$ is possibly with boundary) and with
$p\le n-1$. A {\it Poincar\'e
embedding} of a map  $f\:P\to N$ is
a commutative diagram
$$
\begin{CD}
A @>>> C @<<< \partial N\\
@VVV @VVV\\
P @>>f > N\, ,
\end{CD}
$$
in which
\begin{itemize}
\item The square in the diagram is a homotopy pushout.
\item The map $A \to P$ has the homotopy type of a $(n-p-1)$-spherical
fibration over $P$ (in particular its mapping cylinder $\bar P$
is a Poincar\'e space of dimension $n$ with boundary $A$).
\item The composite $\partial N \to C \to N$ is the inclusion.
\item The image of a fundamental class 
under the composite
$$
H_n(N,\partial N;{\cal L}) \to H_n(N,C;{\cal L}) \cong H_n(P,A; {\cal L}_{|P})
$$
equips $P$ with the structure of a Poincar\'e space with boundary $A$.
Similarly, the image of a fundamental class under
$$
H_n(N,\partial N;{\cal L}) \to H_n(N,P\amalg \partial N;{\cal L}) 
\cong H_n(C,A\amalg \partial N; {\cal L}_{|C})
$$
equips $C$ with the structure of a Poincar\'e space with boundary 
$A\amalg \partial N$.
\end{itemize}

We are not necessarily assuming in the
above that $A\to P$ and $A\amalg \partial N \to C$ are
inclusions. However, we are abusing notation slightly when writing
$H_n(N,C;{\cal L})$; the reader should substitute the appropriate
mapping cylinder in such cases.
\medskip

Roughly, a {\it concordance} of two Poincar\'e embeddings
of $f$, is a map of their associated diagrams which is 
a weak equivalence at each space of the diagram
and which moreover, is the identity at $P,N$ and $\partial N$. We say two
embeddings of $f$ are {\it concordant} if there is a finite chain
of concordances connecting them (see \cite{haef2} for a more detailed
definition).

\subsection*{Spectra} 
A {\it spectrum} $X$ is a
collection of based spaces $\{X_i\}_{i \in {\Bbb N}}$
together with based maps $\Sigma X_i \to X_{i{+}1}$
where $\Sigma X_i$ denotes the reduced suspension of $X_i$.
A {\it map of spectra} $X \to Y$ consists of maps
$X_i \to Y_i$ which are compatible with the structure maps.

Let $G$ be a topological group. For technical reasons,
we assume that the underlying space of $G$ 
is the retract of a cell complex.
A {\it (naive) $G$-spectrum} consists of a spectrum
$X$ such that each $X_i$ is a based (left) $G$-space and 
each structure map $\Sigma X_i \to X_{i{+}1}$ is equivariant,
where it is understood that $G$ acts trivially on the
suspension coordinate of $\Sigma X_i$. A {\it map}
of $G$-spectra  is a map of spectra that
is  compatible with the $G$-action.
A {\it weak equivalence}  of
$G$-spectra is a morphism inducing an isomorphism on
homotopy groups in every degree.
 
S.\ Schwede has shown that $G$-spectra form a Quillen
model category  with the above notion of
weak equivalence\cite[\S2]{Schwede}. 
For reasons of space, we shall
not describe the entire structure,
but we will describe the fibrant and cofibrant
objects arising from Schwede's model structure.
A $G$-spectrum $X$ is {\it fibrant} if it is an $\Omega$-spectrum, meaning
that the adjoint maps $X_j \to \Omega X_{j+1}$ are weak equivalences
for all $j$. It is {\it cofibrant} if
each $X_j$ is built up from a point
by attaching free $G$-cells and each structure map $\Sigma X_j \to X_{j+1}$ is
is obtained up to isomorphism by attaching free cells to $\Sigma X_j$.
We also consider retracts of such an $X$ to be cofibrant.

The {\it suspension spectrum} $\Sigma^{\infty} X$ of a based
$G$-space $X$ has  $j$-th space $Q(S^j \smsh X)$,
where $Q = \Omega^{\infty}\Sigma^{\infty}$ is the stable
homotopy functor. 
\medskip

If $X$ is a  $G$-spectrum then
the {\it homotopy orbit spectrum} $X_{hG}$
is the (non-equivariant) spectrum
given by 
$$
 X \smsh_G (EG_+)\, ,
$$ 
where $EG$ is 
the free contractible $G$-space (arising from the bar construction),
and $EG_+$ is the result of adding a basepoint to $EG$.
This spectrum has $j$-th space $X_j \smsh_G (EG_+)$.

The {\it homotopy fixed point spectrum} $X^{hG}$ is given by
$$
F(EG_+,X)^G \, .
$$
This is the spectrum whose $j$-th space is the equivariant function space
of based maps $EG_+ \to X_j$. Taking homotopy fixed points
is homotopy invariant when $X$ is fibrant.

As in \cite{Klein_dualizing}, we use
naive versions of the smash product (which are
associative, unital and commutative up to homotopy).

\section{The dualizing spectrum of an inclusion}
Let $G$ be a topological group. The {\it group ring} of $G$ over the sphere
spectrum is the $(G\times G)$-spectrum
$$
S^0[G]
$$
which is the suspension spectrum of $G_+$. 
The action of the left copy of $G$ in 
$G\times G$ on $S^0[G]$ is given by left multiplication. 
The action of the right copy is
given by right multiplication composed with inversion ($g\mapsto g^{-1}$).

In \cite{Klein_dualizing}, the {\it dualizing spectrum} of $G$
$$
D_G = S^0[G]^{hG}
$$
was defined as the homotopy fixed points of the left copy of $G$ acting
on the group ring. Since the two $G$-actions commute, the right copy of
$G$ acts on $D_G$ giving it the structure of a 
$G$-spectrum.

If $h\:H \to G$ is a homomorphism, then $S^0[G]$ inherits an
$(H\times G)$-action (restriction of scalars), and we have maps
$$
\begin{CD} 
D_H @>h_! >> S^0[G]^{hH} @< h^! << D_G \, .
\end{CD}
$$
The first of these, called {\it pushforward,} 
is induced by the map $S^0[H] \to S^0[G]$ arising from $h$. 
The second, called {\it restriction} 
is given by regarding a $G$-fixed point as an $H$-fixed point. 
The pushforward map is $H$-equivariant, and
the restriction map is $G$-equivariant.

From this point on,
we will only need consider the situation where
$h\: H \to G$ is an inclusion.

\begin{defn} The {\it dualizing spectrum} associated
with an inclusion $H\subset G$ is the $G$-spectrum
$$
D_{H\subset G} \,\, := \,\, S^0[G]^{hH} \, .
$$
\end{defn}

Before proceeding further, it will be useful to recall some results
from \cite{Klein_dualizing}.

\begin{thm}[{\cite[cor. 5.1]{Klein_dualizing}}] 
\label{dual} Assume $BG$ is a finitely dominated space. 
Then
\begin{itemize}
\item $D_G$ is a suspension spectrum, i.e., 
there are an integer $j\ge 0$, a finitely dominated, $1$-connected
based $G$-space $Y$ and an equivariant weak equivalence
$$
\Sigma^j D_G \simeq \Sigma^\infty Y \, .
$$
\item If furthermore $BG$ is a Poincar\'e space (of dimension $d$), then $Y$ is
unequivariantly weak equivalent to a sphere (of dimension $j-d$), 
and the Spivak normal fibration of $BG$ is given by the Borel construction
$$
Y \to EG\times_G Y \to BG \, .
$$
\end{itemize}
\end{thm}

We need a partial extension of \ref{dual} to $D_{H\subset G}$.

\begin{prop} \label{susp} Assume $BH$ is finitely dominated.
Then the first conclusion of \ref{dual} holds for $D_{H\subset G}$
(i.e., it is a suspension spectrum).
\end{prop}

\begin{proof} 
The proof of the proposition is basically the
same as that of the first part of \ref{dual}.
We will merely sketch the argument.

Since $BH$ is finitely dominated,
the homotopy quotient $EH\times_H G$  is $G$-finitely dominated (since the
Borel construction of $G$ acting on it is identified with $BH$).
By the equivariant duality theory developed in \cite{Klein_immersion}, there
is an integer $j\ge 0$,  
a 1-connected finitely dominated based free $G$-CW complex $Y$ and a map
$$
d\:S^j \to Y \smsh_G (EH\times_H G)_+
$$
which satisfies duality, meaning that for all cofibrant and fibrant
$G$-spectra  $E$, the map
$d$ induces an equivalence of mapping spectra
$$
\text{map}((EH\times_H G)_+,E)^G \,\, \simeq\,\,  
\text{map}(S^j,Y\smsh_G E) \, ,
$$ 
In particular, if $E$ has a 
$(G\times G')$-action, the above equivalence is $G'$-equivariant.

By `change of groups,' the domain of this equivalence is
identified with $E^{hH}$.
Setting $E= S^0[G]$ (considered as a $(G\times G)$-spectrum) the equivalence 
of mapping spectra therefore becomes
$$
D_{H\subset G} \,\, \simeq \,\, S^{-j}\smsh Y \, ,
$$
and it is $G$-equivariant.
Applying $j$-fold suspension to both sides completes the proof.
\end{proof}

\section{Proof of Theorem \ref{embed_thm} when $P$ is connected}
A choice of  basepoint for
$P$ makes $f\: P \to N$ a map of based spaces. 
If $P$ is connected, then $f(P) \subset N$
is contained in a connected component. It will be enough
to produce a Poincar\'e embedding into that connected component.
For this reason, we can assume without loss in generality that
$N$ is also connected.
\medskip

Let $H$ be a topological
group model for the loop space $\Omega P$ and  $G$ be one for $\Omega N$.
Here is one way to get such models: let $S_\cdot N$ denote the simplicial
total singular complex of $N$, and let $G_\cdot$ denote its Kan loop group.
Define $G$ to be the geometric realization of the underlying simplicial set
of $G_{\cdot}$. Then $N \simeq BG$. Furthermore, the construction is natural,
so we get a homomorphism $H \to G$ whose classifying map 
$BH \to BG$ is identified with $f$. 

As it stands, the homomorphism $H\to G$ isn't necessarily an inclusion.
However, one can always achieve this by replacing the map 
$f\: P \to N$ by its  mapping cylinder.  
Assume this has been arranged.
\medskip

As in the previous section, there are pushforward and restriction maps
$$
D_H \to D_{H\subset G} \leftarrow D_G\, ,
$$
where
the first of these is $H$-equivariant and the second is $G$-equivariant.
Applying \ref{susp}, there is an integer $j \ge 0$, a $1$-connected finitely
dominated based $Y_H$ with $H$-action, $1$-connected finitely
dominated based $G$-spaces $Y_{H\subset G}$  and $Y_G$ and equivariant
weak equivalences
$$
\Sigma^\infty Y_H \simeq \Sigma^j D_H\, , \quad 
\Sigma^\infty Y_{H\subset G} \simeq \Sigma^j D_{H\subset G}\, 
\quad \text{ and } \quad 
\Sigma^\infty Y_G \simeq \Sigma^j D_G\, . 
$$
Since  $BH$ is a closed Poincar\'e space, $Y_H$ has the
unequivariant homotopy type of a sphere.

With respect to these identifications, the pushforward and restriction
maps, after $j$-fold suspending, are identified with {\it stable} maps
$$
Y_H \to Y_{G\subset H} \leftarrow Y_G \, ,
$$
the left of these is $H$-equivariant and the right one is $G$-equivariant.
We may assume that each of the three spaces above is a 
equivariant based CW complex which is free away from the basepoint.
Since all of these are equivariantly finitely dominated, 
a straightforward
obstruction theory argument shows that, at the expense of suspending
these spaces finitely many times, we can represent the above stable
maps by unstable ones (i.e., we are in effect replacing $Y_H,Y_{H\subset G}$
and $Y_G$ by their finite $k$-fold suspensions for some $k$
sufficiently large.) Therefore, we
can assume without loss in generality that the above equivariant stable maps 
are {\it unstable.}

Applying the Borel constructions,
we obtain a commutative diagram
\begin{equation}\label{embedding_diagram}
\xymatrix{
E_H \ar[r]\ar[d] & E_{H\subset G} \ar[d] & 
E_G \ar[l]  \\
BH \ar[r] & BG
}
\end{equation}
in which $E_H = Y_H \times_H EH$, etc. 

Let  $T_H$ and $T_G$ be the mapping cylinders
of the maps $E_H \to BH$ and $EG \to T_G$.
Then
$$
(T_H,E_H) \quad \text{ and } (T_G,E_G)
$$
are Poincar\'e pairs of the same dimension 
(see \cite[p.\ 22]{Klein_dualizing}).

\begin{prop}\label{embed_framed} With respect to the identification
$T_G \simeq BG$, the diagram  above is a Poincar\'e embedding
of the map
$$
BH \to BG \overset\sim \to T_G
$$
\end{prop}

\begin{proof} We need to first show that the square
\begin{equation}\label{pushout}
\xymatrix{
E_H \ar[r]\ar[d] & E_{H\subset G} \ar[d]   \\
BH \ar[r] & BG
}
\end{equation}
is a homotopy pushout. 

To see this, first note that the vertical fibrations of the square
come equipped with sections, since
they are unreduced Borel constructions of the based spaces $Y_H$ and
$Y_{H\subset G}$. The homotopy
pushout property will follow if we can show that the induced map between
the quotient spaces induced by these sections is a weak equivalence.

Now, if $Y$ is a based $G$ space, then the quotient space
of the canonical section $BG \to Y \times_G EG$ is just the {\it reduced}
Borel construction
$$
Y\smsh_G (EG)_+\, .
$$
Consequently, we need to show that the map of reduced 
Borel constructions
$$
Y_H \smsh_H (EH)_+ \to Y_{H\subset G}\smsh_G (EG)_+
$$
is a weak equivalence. This is a map of $1$-connected spaces,
so it will be enough to check that the induced map of suspension
spectra is a weak equivalence. After desuspending, 
the map of suspension spectra
is identified with the map of homotopy orbit spectra
$$
(D_H)_{hH}\to (D_{H\subset G})_{hG} \, .
$$
In order to show that the latter map is an equivalence, we  
use the {\it norm map}
$$
\eta_H \: D_H \smsh_{hH} E \to E^{hH}
$$
This map was introduced
in \cite{Klein_dualizing} and was shown to be a weak equivalence for all
$E$ whenever $BH$ is a finitely dominated space. The norm map
is natural in $E$.

Associated with an inclusion $H\subset G$ there is, more generally, 
a norm map
$$
\eta_{H\subset G} \:D_{H\subset G} \smsh_{hG} E \to E^{hH}
$$
which is natural in $G$-spectra $E$ (the construction of this
more general norm 
map is virtually identical to the one given in 
\cite{Klein_dualizing} and we will leave these details to
the reader). Again this norm map is a weak equivalence whenever 
$BH$ is finitely dominated.

Furthermore, both kinds of norm maps
are compatible with pushforward and restriction in
the sense that the evident diagram
\begin{equation}\label{norm}
\xymatrix{
D_H\smsh_{hH} E \ar[rr]^{\eta_H}_\sim  \ar[d] &&  E^{hH} \ar^{||}[d]\\
D_{H\subset G} \smsh_{hG} E \ar[rr]^{\eta_{H\subset G}}_\sim && E^{hH} \\
D_G\smsh_{hG} E \ar[rr]_{\eta_G}^\sim\ar[u] && E^{hG} \ar[u]  
}
\end{equation}
is homotopy commutative. Applying the upper 
square of the diagram to $E = S^0$
(with trivial action) shows that the map
$$
(D_H)_{hH} \to (D_{H\subset G})_{hG}
$$
induced by pushforward is a weak equivalence. We have now shown that
\eqref{pushout} is a homotopy pushout.

The second step of the proof is to establish compatibility of the fundamental
classes. By \cite[2.3]{haef}, we are only required to show that the fundamental
classes of $(T_G,E_G)$ and $(T_H,E_H)$ are compatible. 
For this, we need to know how the fundamental classes arise.

Recall that  $E_G \subset T_G$ is stably identified with
$D_G \times_G EG \to BG$. Up to a shift in degree, the homology
of the quotient $T_G/E_G$ is the same as the homology of the spectrum
given by collapsing the preferred section $BG \to E_G\times_G EG$ and
this is just the homotopy orbit spectrum $(D_G)_{hG}$. Let
$$
\alpha_G\: S^0 \to (D_G)_{hG}
$$
be the homotopy class which corresponds to the `counit' map 
$$BG_+ \to S^0$$ 
(given unstably  by mapping $BG$ to the non-basepoint
of $S^0$)  
by means of the norm equivalence
$$
\eta_G\: (D_G)_{hG} \simeq D_G \smsh_{hG} S^0 \overset \sim \to 
(S^0)^{hG} = \text{map}(BG_+,S^0)\, .
$$
Then with respect to our identifications,
$\alpha_*([S^0]) \in H_0((D_G)_{hG})$ is a fundamental class
(cf.\  \cite[p.\ 441]{Klein_dualizing}).

Compatibility of the fundamental classes is now a direct consequence
of the homotopy commutativity of \eqref{norm} and the observation that
the composite $BH_+\to BG_+ \to S^0$ is the counit map.
The completes the proof of \ref{embed_framed}.\end{proof}

\subsection*{Untwisting}
By Proposition \ref{embed_framed} 
$$
\begin{CD}
BH @>>> BG @> \sim> > T_G
\end{CD}
$$
Poincar\'e embeds with diagram
\eqref{embedding_diagram}. Although $T_G$ has the homotopy type of $N$,
the Poincar\'e pair $(T_G,E_G)$ need not have the homotopy type 
of $(N\times D^j,\partial (N \times D^j))$. To obtain a Poincar\'e embedding
into the latter, we will need to twist by an inverse to the Spivak fibration
of $N$.

Let $\tau\:S(\tau)\to N$ be a choice of fiber homotopy
inverse to the Spivak fibration of $N$. Continuing to identify
$N$ with $BG$, we form the fiberwise join
$$
S(\tau) \ast_{BG} E_G \to BG \, .
$$
To simplify notation denote its total space by
${}^\tau E_G$.

Similarly, applying fiberwise join with $\tau$ to the
other terms in the
top row of \ref{embedding_diagram} we obtain a modified diagram
$$
\xymatrix{
{}^\tau E_H \ar[r]\ar[d] & {}^\tau E_{H\subset G} \ar[d] & 
{}^\tau E_G \ar[l]  \\
BH \ar[r] & BG \, .
}
$$
which represents a Poincar\'e embedding, but this time
the ambient Poincar\'e space is $N \times D^j$.
The proof of Theorem \ref{embed_thm} is now complete in the case when
$P$ and $N$ are connected.

\section{The fibered theory}
Recall that $G$ is a topological group model for the loop
space $\Omega N$.
The unreduced Borel construction of $G$ acting on
the dualizing spectrum $D_G$ produces a fibered spectrum
over $BG$.
A more direct construction can be given producing
a fibered spectrum over $N$. 
The latter will enable us to complete 
the proof of Theorem \ref{embed_thm} in
the general case. The new construction is
conceptually simpler, albeit more technical.
\medskip

\begin{defn} Let $X$ be a space. Let  $R(X)$
be the category of {\it retractive spaces} over $X$. An {\it object} 
is a space $Y$ equipped with maps $s_Y\:X\to Y$ and $r_Y\: Y\to X$
such that $r_Y\circ s_Y$ is the identity map. 
A {\it morphism} $f\:Y\to Z$ is a
map of underlying spaces which commutes with their given structure maps:
$r_Z\circ f = r_Y$ and $f\circ s_Y = s_Z$.
A morphism is a {\it weak equivalence} 
if it is a weak homotopy equivalence of underlying spaces.
It is a {\it fibration} if it is a Serre fibration of underlying spaces.
It is a {\it cofibration} if it has the left lifting property with 
the acyclic fibrations (= fibrations which are also weak equivalences). 
In particular, an object $Y$ is {\it fibrant} if $r_Y\: Y\to X$ is a Serre fibration. It is {\it cofibrant} 
if $s_Y\: X \to Y$ is a Serre cofibration (i.e., 
$s_Y$ is an inclusion and the pair $(Y,X)$ is a retract
of a pair $(Z,X)$ in which $Z$ is obtained from $X$ by attaching
cells).
\end{defn}

\subsection*{Fibered spectra}
S.\ Schwede has shown that the spectra built from objects of
any (pointed, simplicial) model category again form a model category.
We will apply his result to the model category of retractive spaces
$R(X)$. 

Given objects $Y,Z\in R(X)$, the hom-set $\hom_{R(X)}(Y,Z)$ may be topologized as a subspace
of the function space of all continuous maps $Y \to Z$ of underlying spaces,
where the function space is equipped with the compactly generated compact
open topology. This gives 
$R(X)$ the structure of a {\it topological} model category.
In particular, $R(X)$ is simplicial.
\medskip

Given an object $Y \in R(X)$, its 
{\it (reduced) fiberwise suspension} is given by 
$$
\Sigma_X Y \,\, = \,\, (Y \times D^1) \cup_{(X\times D^1)\cup (Y\times S^0)} X\, ,
$$
where $(X\times D^1)\cup (Y\times S^0)$ is amalgamated
along $X\times S^0$ and the map 
$(X\times D^1)\cup (Y\times S^0) \to X$
 is defined as the composite of the inclusion
$(X\times D^1)\cup (Y\times S^0) \subset Y\times D^1$,
followed by first factor projection $Y \times D^1 \to Y$, followed
by $r_Y \: Y \to X$.
Then  $\Sigma_X$ defines an endo-functor of $R(X)$.

\begin{defn} A {\it fibered spectrum} 
$ {\cal E}$ over $X$
consists of  objects ${\cal E}_j \in R(X)$ for $j \in {\Bbb N}$ together
with (structure) maps
$$
\Sigma_X {\cal E}_j \to {\cal E}_{j+1} \,  ,
$$
where $\Sigma_X\: R(X) \to R(X)$ is the reduced fiberwise suspension functor.
A {\it morphism} ${\cal E}\to {\cal E}'$ is given by   
maps ${\cal E}_j \to {\cal E}'_j$ which are compatible with the
structure maps.

We say that ${\cal E}$ is {\it fibrant} if the adjoints to the structure
maps are weak homotopy equivalences of underlying spaces.
Any fibered spectrum ${\cal E}$ can be converted into a fibrant
one ${\cal E}^{\text{f}}$ in which 
$$
{\cal E}^{\text{f}}_j\,\, := \,\, \underset{j}{\text{hocolim\, }} 
\Omega^j_X {\cal E}_j \, ,
$$
where the homotopy colimit is taken in $R(X)$, and 
$\Omega_X^j$ is the adjoint to $j$-fold reduced fiberwise
suspension. The above is called {\it fibrant 
replacement.}

We are now ready to describe the model structure on fibered spectra.
A morphism ${\cal E}\to {\cal E}'$ is a 
{\it weak equivalence} if the associated morphism of 
fibrant replacements  ${\cal E}^\text{f}\to ({\cal E}')^\text{f}$ 
is a {\it levelwise} weak equivalence: 
for each $j$, the map ${\cal E}_j^\text{f}\to ({\cal E}')_j^\text{f}$ 
is required to be weak equivalence of $R(X)$. 
A morphism ${\cal E}\to {\cal E}'$ is a 
{\it cofibration} if the maps 
$$
{\cal E}_0 \to {\cal E}'_0 \quad \text{ and } \qquad
{\cal E}_j \cup_{\Sigma_X {\cal E}_{j-1}}  \Sigma_X {\cal E}'_{j-1}\to
 {\cal E}'_j 
$$
are cofibrations of $R(X)$. A morphism is a {\it fibration} if it has the right
lifting property with respect to the acyclic cofibrations.
\end{defn}

\begin{rem}  
Several model structures on parametrized spectra are described
in the forthcoming book by May and Sigurdsson
\cite{May-Sigurdsson}. This is one of them.
\end{rem}

\subsection*{The fibered spectrum ${\cal D}(X)$}
Let $(X,x)$ be a based space and $j \ge 0$ an integer. 
Then the wedge $X\vee_x S^j$ is defined and we have a projection
$$
X \vee_x S^j \to X \, .
$$
We convert the latter into a Hurewicz fibration by replacing the
wedge with the homotopy equivalent space
$$
E(X,x)_j := X^I \cup_{P_x X} P_x X \times S^j \, , 
$$
where $X^I$ is the space of paths $[0,1]\to X$ and $P_x X$ is the subspace
of those paths which map $0$ to $x$. The inclusion $P_x X \to P_x X \times S^j
$ is given by selecting the basepoint of $S^j$.
The map $E(X,x)_j \to X$ is induced by evaluating a path at its endpoint.
Note that $E(X,x)_j \to X$ has a preferred section given by
$x \mapsto c_x$, where $c_x$ is the constant path having value $x$.

\begin{lem} The fiber at $x$ of $E(X,x)_j \to X$ 
is, up to homotopy,
the smash product
$$
(\Omega_x X)_+ \smsh S^j \, ,
$$
where $\Omega_x X$ is the based loop space of $X$ and 
$(\Omega_x X)_+$ is the effect of adding a disjoint basepoint.
\end{lem}

\begin{proof} The fiber at $x\in X$ is 
$$
P_x X \cup_{\Omega_x X} (\Omega_xX \times S^j)
$$
and this has the homotopy type of the quotient
$$
(\Omega_xX \times S^j)/(\Omega_xX \times *) \,\, = \,\, 
(\Omega_xX)_+ \smsh S^j\, .
$$
\end{proof}

\begin{defn} Let $${\cal D}(X,x)_j$$ be the space of sections
of $E(X,x)_j \to X$. Note that  ${\cal D}(X,x)_j$  comes equipped with
a preferred basepoint. 

Define a fibration 
$$
{\cal D}(X)_j \to X
$$
whose total space is the space of pairs 
$$
(x,\sigma)
$$
in which $x \in X$ and $\sigma\in D(X,x)_j$ (topologize
it as a subspace of $X\times \text{\rm map}(X,X^I \times S^j)$).
As each fiber ${\cal D}(X,x)_j$ is based,
${\cal D}(X)_j \to X$ comes equipped with a preferred section.
\end{defn}

The collection 
$$
{\cal D}(X) \,\, := \,\, \{{\cal D}(X)_j\}_j
$$
can be given the structure of a fibered spectrum over $X$. To see this,
define an auxiliary space ${\cal C}(X)_j$ in a way similar to 
${\cal D}(X)_j$, but where we replace the sphere $S^j$ by $D^{j+1}$.
The standard radial contraction of the disk to a point induces
a fiberwise contraction of ${\cal C}(X)_j$ to the zero object.

Furthermore,
there is an evident pushout square of inclusions
$$
\begin{CD}
{\cal D}(X)_j @>\subset >> {\cal C}(X)_j\\
@V\cap VV @VV\cap V \\
{\cal C}(X)_j @>>\subset > {\cal D}(X)_{j+1} \, .
\end{CD}
$$
Using the null homotopies provided by 
each copy of ${\cal C}(X)_j$, we obtain
a map $\Sigma_X {\cal D}(X)_j \to  {\cal D}(X)_{j+1}$.
We infer that that ${\cal D}(X)$
is a fibered spectrum.

The following result explains the relationship between
${\cal D}(X)$ and the dualizing spectrum of the loop space
of $X$. 

\begin{thm} If $X \simeq BG$, then there is a weak equivalence of
fibered spectra
$$
{\cal D}(X) \,\, \simeq \,\, EG\times_G D_G\, .
$$
\end{thm} 

\begin{proof} The proof will use the following fact: up to fiber homotopy
equivalence, a fibration over a connected based space can always
be regarded as a Borel construction. Specifically, if $E\to X$ is a
fibration, we identify $X$ with $BG$ for a suitable topological group $G$.
The fiber product $$EG \times^{BG} E$$ is then a $G$-space (because it is
a $G$-invariant subspace of the product $EG \times E$, where
$G$ acts trivially on $E$) whose
Borel construction recovers $E$ up to fiberwise weak
equivalence. Notice that $EG \times^{BG} E$
has the underlying homotopy type of the fiber of $E\to X$;
we call it the {\it thick fiber}.

We now begin the proof.
Without loss in generality, take $X = BG$.
Consider first the fibration
$$
{\cal D}(BG)_0\to BG\, .
$$
Let $x\in BG$ be a point.
Then the fiber at $x$ is the space of sections
of the fibration
$$
BG^I \cup_{P_xBG} (P_xBG) \times S^0 \to BG\, .
$$
The latter  has a preferred identification 
up to homotopy with the fibration
$$
BG \amalg E_xG \to BG \, .
$$
where $E_x G = EG\times_G G_x$, and $G_x$ is the
fiber of the universal principal bundle $EG \to BG$ at
$x$. 

The space of sections of this last fibration is precisely
the homotopy fixed point space
$$
((G_x)_+)^{hG} \,\, = \,\, \text{maps}(EG,(G_x)_+)^G\, .
$$

We now compute the thick fiber. 
The above shows that it is given up to equivariant
weak equivalence by
the space of pairs
$$
(y,\sigma)
$$
with $y \in EG$ and $\sigma$ a homotopy fixed point of
$$
((G_x)_+)^{hG}\, ,
$$
where $x \in BG$ is the projection of $y$ (i.e., $y \in G_x$). 
By definition of the thick fiber,
the action of $G$ on this pair is given by $g\cdot (y,\sigma) :=
(gy,\sigma)$.

Since $y \in G_x$, we have an isomorphism $$h_y\:G \overset\cong\to G_x$$ 
given by $g\mapsto gy$. Using this isomorphism, we
get a homotopy fixed point 
$$
\sigma_y = (h_y)^{-1}\circ \sigma\in (G_+)^{hG} \, .
$$ 
The assignment
$$
(y,\sigma) \mapsto \sigma_y
$$
defines an equivariant weak equivalence from the thick fiber
to $(G_+)^{hG}$. 
We infer that ${\cal D}(BG)_0$ is fiber homotopy
equivalent to  
the Borel construction of $G$ acting on $(G_+)^{hG}$.

The case of the fibration 
$$
{\cal D}(BG)_j \to BG
$$
for $j > 0$ is essentially the same, but now the
argument gives an equivariant weak equivalence between
the thick fiber and $(S^j \smsh (G_+))^{hG}$. 
Hence, we have
a fiber homotopy equivalence
$$
{\cal D}(BG)_j \,\, \simeq \,\, EG \times_G (S^j \smsh (G_+))^{hG} \, .
$$
We leave it to the reader to verify that the equivalences for each $j$
that we produced give a morphism of fibered spectra.
We have now shown that the fibered spectrum
${\cal D}(BG)$ is identified with the fibered spectrum
whose $j$-th term is $EG \times_G (S^j \smsh (G_+))^{hG}$. 
To complete the proof, we identify this last fibered spectrum 
with the unreduced Borel
construction of $D_G$. 

This follows immediately from the fact that
the evident
map of spaces
$$
EG\times_G (S^j \smsh (G_+))^{hG} \to 
EG \times_G Q(S^j \smsh (G_+))^{hG}
$$
is $(2j-c)$-connected, for some constant $c\ge 0$ independent of $j$
(this uses the Freudenthal suspension theorem
and the assumption of  $X = BG$ being finitely dominated).
A map of such fibered spectra is
clearly a weak equivalence.
\end{proof}

\begin{rem} If $X = \amalg_\alpha  X_\alpha$ 
is a space decomposed into connected components, then
we have a decomposition of the fibered spectrum
$$
{\cal D}(X) \,\, = \,\, \coprod_\alpha {\cal D}(X_{\alpha})\, .
$$
In particular, if $X$ is a Poincar\'e space, the above theorem shows that
${\cal D}(X) \to X$ is the stable Spivak fibration.
\end{rem}

\subsection*{Extension to maps} Given a map of spaces
$$
A \to X \, 
$$
we can associate a fibration
$$
{\cal D}(A\to X)_j \to X
$$
whose total space consists of pairs
$$
(x,\sigma)
$$
in which $x\in X$ an $\sigma$ is a section of $E(X,x)_j \to X$ {\it along} $A$,
i.e., $\sigma$ is a section of the associated pullback fibration.
The collection
$$
{\cal D}(A\to X)\,\,  := \,\, \{{\cal D}(A\to X)_j\}_{j\ge 0} \, ,
$$
forms a fibered spectrum over $X$.
The case of the identity map $X\to X$ recovers ${\cal D}(X)$.

If $f\:A\to B$ is a map of spaces
over $X$, then we obtain a {\it restriction} map
$$
f^! \: {\cal D}(B\to X) \to {\cal D}(A\to X)
$$
This is a morphism of fibered spectra over $X$.

There is also a {\it pushforward} map
$$
f_!\: {\cal D}(A\to B) \to {\cal D}(A\to X)
$$
which is a map of fibered spectra covering $B\to X$.

\section{Proof of Theorem \ref{embed_thm} in the general case}
We apply the constructions of the 
last section to $f\: P \to N$ to get a diagram
$$
\xymatrix{
{\cal D}(P) \ar[r]^{f_!\quad}\ar[d] & {\cal D}(P\to N) \ar[d] & 
{\cal D}(N)\ar[l]_{\quad f^!} \ar[d] \\
P \ar[r]_f & N & N\ar^{=}[l]& \, .
}
$$
By analogy with the connected case, each of the fibered spectra
appearing above is fibered suspension spectrum
(we omit the details).

The rest of the proof is as in the connected case:
the above diagram 
fiberwise suspends to an associated commutative diagram of spaces
$$
\xymatrix{
E_P \ar[r]\ar[d]  & E_{P \to N}  \ar[d] & E_N \ar[l]\\
P \ar[r]_f & N\, , 
}
$$
which 
is a Poincar\'e embedding into the mapping cylinder
of ${E_N \to N}$. We must then twist if necessary by a
suitable element 
in the Grothendieck group of spherical fibrations
over $N$ to get a Poincar\'e embedding into $N \times D^j$
for suitable $j$. The proof of Theorem \ref{embed_thm} is
now complete. 

\section{Proof of Theorem \ref{uniqueness}}
We will prove uniqueness when $N$ is closed
and has trivial Spivak fibration. The general case, which follows
by a twisting argument similar to that of the last section,
is left to the reader.

Suppose we are given two 
Poincar\'e embeddings of $f_j\: P \to N\times D^j$, with 
diagrams
$$
\xymatrix{
A_i \ar[r]\ar[d] & C_i \ar[d]& N \times S^{j-1} \ar[l]\\
P \ar[r]_{f_j\quad } & N\times D^j 
}
\qquad i = 0,1 \, .
$$
the assumptions on $N$ imply that $A_i \to P$ are Spivak fibrations for $P$.

To minimize technical difficulties, it is more convenient to replace 
$N \times D^j$ by $N$ (by means of first factor projection), 
thereby rewriting the diagrams as
$$
\xymatrix{
A_i \ar[r]\ar[d] & C_i \ar[d]& N \times S^{j-1} \ar[l]\\
P \ar[r]_{f} & N
}
\qquad i = 0,1 \, .
$$
There is no loss of information, since
we can recover the original diagrams up to concordance
by taking suitable mapping cylinders.

If $j$ is large, the Spivak fibration is unique up to fiber homotopy
equivalence. Hence there is a fiber homotopy equivalence over $P$
$$
A_0 \simeq A_1 \, .
$$
Therefore, without loss in generality, we can assume $A_0 = A_1$.
To simplify notation, we set 
$$
A \,\, := \,\, A_0 \, .
$$
The next step is to relate $C_0$ with $C_1$. 
We first make some definitions.

\begin{defn} Let $T(N)$ be the category whose {\it objects}
are spaces $Y$ equipped with a map $Y \to N$. A {\it morphism}
is a map compatible with the structure map to $N$.
Thus $T(N)$ is the category of spaces over $N$.
A morphism is a {\it weak equivalence} when the underlying
map of spaces is a weak homotopy equivalence.

The {\it unreduced fiberwise suspension} of $Y$ is defined
to be
$$
S_N Y \,\, := \,\, (Y\times D^1) \cup_{Y\times S^0} (N \times S^0) \, .
$$
It is an object of $T(N\times D^1)$ but we usually regard it as
an object of $T(N)$ using the first factor projection $N\times D^1 \to N$.
\end{defn}

\begin{defn}  Given a Poincar\'e embedding
$$
\xymatrix{
A \ar[r]\ar[d] & C \ar[d] & N\times S^{j-1} \ar[l]\\
P    \ar[r]_{f} & N\, ,
}
$$
its {\it decompression} is the associated Poincar\'e embedding 
$$
\xymatrix{
\Sigma_P A \ar[r]\ar[d] & S_N C \ar[d] & N\times S^j \ar[l]\\
P  \ar[r]_{f} & N \, ,
}
$$
where  $N\times S^j \to S_N C$ is the fiberwise
suspension of $N \times S^{j-1} \to C$.
\end{defn}

\begin{defn}
Let $A\to X$ be an inclusion map of  $T(N)$. 
The {\it fiberwise 
quotient} is
$$
X/\!\!/A := X \cup_A N \, , 
$$
which is  an object of the retractive space category $R(N)$, but can also
be regarded as an object of $T(N)$ by means of the  forgetful
functor $R(N) \to T(N)$.
\end{defn}

To see that $S_N C_0$ and $S_N C_1$ are
weak equivalent, we compare them both to $P/\!\! /A$ by means of
the `excision' weak equivalences
$$
\begin{CD}
P/\!\! /A := P \cup_A N  @>\sim >> N \cup_{C_i} N \simeq S_N C_i \, .
\end{CD}
$$
This produces the desired chain of weak equivalences.

The final step in the proof of uniqueness is to consider the
maps
$$
N\times S^j \to S_N C_i
$$
appearing in the decompressions of our given Poincar\'e embeddings. 
Call these maps $\alpha_i$.
We must show that our identification $S_N C_0 \simeq S_N C_1$ 
is compatible with $\alpha_i$ up to homotopy.

This identification arises by identifying
each $S_N C_i$ with $P/\!\!/A$. So the $\alpha_i$
correspond to a pair of maps
$$
\beta_i\: N\times S^j \to P/\!\!/A
$$
which we need to prove are fiberwise homotopic. An unraveling of the
definitions shows that these are the fiberwise {\it collapse
maps} of  the two Poincar\'e embeddings (see \S8).

Each one of these maps is therefore {\it fiberwise dual} to the
map
$$
f^+\: P^+ \to N^+
$$
(again, see \S8).
By the uniqueness of fiberwise duality maps, it follows
that $\beta_0$ and $\beta_1$ are homotopic when $j \gg 0$.
This completes the proof of uniqueness.

\section{Appendix: fiberwise duality}

Suppose $N$ is a Poincar\'e duality space, possibly with
boundary $\partial N$. Let $X$ and $Y$ be objects of $R(N)$.
The {\it fiberwise smash product} $X\smsh_N Y$ of $X$ and $Y$
is the object
$$
X \times_N Y \cup_{X \cup_N Y} N 
$$
where $X\times_N Y$ is the fiber product of $X$ and $Y$ 
along $N$. The map $X \cup_N Y \to X\times_N Y$ is
defined using the structure maps $N\to X$ and $N \to Y$.

\begin{rem} As usual, in order to get a homotopy invariant theory,
we must replace the fiberwise smash product with its
derived version, i.e., we replace $X$ and
$Y$ by their fibrant/cofibrant approximations. In what
follows, we will be suppressing this aspect from the notation.
The reader is forewarned.
\end{rem}

Let
$$
N^+ \,\, := N/\!\!/\partial N
$$
denote the {\it double} of $N$. This is 
an object of $R(N)$ whose underlying space is  
$N\cup_{\partial N} N$
Here are some special cases:

\begin{enumerate}
\item If $\partial N$ is empty, then $N^+$ is just $N \times S^0$.
\item The $j$-fold unreduced fiberwise suspension $S^j_N N^+$ coincides with
$(N\times D^j)^+$. In particular, if $\partial N$ is empty, 
we get $N\times S^j$. 
\end{enumerate}

A {\it duality map} for $X$ and $Y$ is a morphism
$$
d\:N^+ \to X\smsh_N Y
$$
Such that for all (cofibrant and fibrant) fibered spectra ${\cal E}$, the operation
$g \mapsto (g\smsh_N {\text{id}_Y})\circ d$ induces an isomorphism
of abelian groups
$$
[X,{\cal E}]_N \cong [N^+,{\cal E}\smsh_N Y]_N
$$
where $[\quad,\quad]_N$ means fiberwise homotopy classes, and
the fiberwise smash product ${\cal E}\smsh_N Y$ 
is defined by $({\cal E}\smsh_N Y)_j = {\cal E}_j \smsh_N Y$.

\begin{rems} 
{\flushleft (1).} When $N= D^n$ is a disk, fiberwise duality
amounts to ordinary Spanier-Whitehead $n$-duality.
\medskip

(2). When $N$ is an arbitrary space, 
a closely related type of fiberwise duality
theory was considered by Vogell \cite{Vogell}, who considered
maps of $R(N\times N)$ of the form
$$
X \,\, \sharp \,\,  Y \to T_j
$$
where `$\sharp$' 
denotes the {\it external} fiberwise smash product of 
$X$ and $Y$ (this is retractive over $N\times N$), and $T_j$
is given as follows: consider $N$ as a space over $N\times N$ via the diagonal
and make it retractive by adding a disjoint copy of $N\times N$. 
Call the resulting object $T_0$.
Then we take $T_j$ to be the 
$j$-fold fiberwise suspension of $T_0$ (in Vogell's terminology,
$T_j$ is the fiberwise Thom complex of the rank $j$-trivial bundle). 
If $N$ is a closed Poincar\'e space, then 
Vogell's duality maps are related to Richter's via a certain 
{\it umkehr} correspondence. The theories in this case are equivalent.
\end{rems}
\medskip

We now state without proof some basic properties of fiberwise duality maps.

\begin{itemize}
\item ({\it Suspension}). If $d$ is a fiberwise duality map, so is its fiberwise suspension
$$
\Sigma_N d\: (N\times D^1)^+ \to X\smsh_N \Sigma_N Y \, .
$$
\item  ({\it Switching}). If $d$ is a fiberwise duality map, so is the composite
$$
\begin{CD}
 N^+ @> d >> X\smsh_N  Y @> \text{twist} >>  Y\smsh_N  X \, .
\end{CD}
$$
\item ({\it Existence}). Given a finitely dominated object $X$, there is always
a finitely dominated object $Y$, an integer $j\ge 0$ and 
a duality map 
$$
(N\times D^j)^+ \to X\smsh_N Y \, .
$$
\item ({\it Uniqueness}). 
Without loss in generality take $Y$ to be both fibrant
and cofibrant. Let
$d'\: N^+ \to X\smsh_N Z$ be another fiberwise duality
map
with $Z$ cofibrant and fibrant. Assume $j$ is large. 
Then there is a weak equivalence
$
h\:\Sigma^j_N Y \overset\sim\to \Sigma^j_N Z
$
such that the composite
$$
\begin{CD}
(N\times D^j)^+ @> \Sigma^j_N d >> X\smsh_N \Sigma^j_N Y @> {\text{id}\smsh_N h} >> X \smsh_N \Sigma^j_N Z
\end{CD}
$$
has the same fiberwise homotopy class as  $S^j_N d'$.
\end{itemize}

Geometrically, fiberwise duality maps arise from Poincar\'e embeddings
in the following way: suppose we are given a Poincar\'e embedding
$$
\xymatrix{
A \ar[r] \ar[d] & C \ar[d] & \partial N \ar[l]  \\
P \ar[r]_{f_j} & N \, .
}
$$
Then there is a  fiberwise duality map
$$
d\:N^+ \to P^+ \smsh_N P/\!\! /A 
$$
which may be defined as a composition of two maps:  the 
fiberwise {\it collapse}
$$
N^+ \to P/\!\! /A
$$
followed by the fiberwise {\it diagonal}
$$
P/\!\! /A \to P^+ \smsh_N P/\!\! /A \, ,
$$
(where $P^+$ denotes $P \amalg N$).

As a fiberwise homotopy class, the 
collapse arises from the chain of fiberwise pairs
$$
\begin{CD}
(N,\partial N) @< \sim << (P \cup_A C,\partial N) @>>> (P \cup_A N,N)
\end{CD}
$$
by taking fiberwise quotients. Similarly,
the fiberwise diagonal  arises from the evident fiberwise
diagonal map of pairs
$$
\begin{CD}
(P,A) @>>> (P \times_N P,P \times_N A ) 
\end{CD}
$$
by taking fiberwise quotients.
\medskip

The following result is due to Richter (unpublished; see
\cite[p.\ 163]{Vogell} for a dual version in the manifold case).

\begin{prop}[Richter Duality]
The map $$d\:N^+ \to P^+ \smsh_N P/\!\! /A$$
is a fiberwise duality.
\end{prop}

Assume now that $\partial N = \emptyset$
and we are given a Poincar\'e embedding
$$\xymatrix{
A \ar[r] \ar[d] & C \ar[d] & \ar[l] N\times S^{j-1} \\
P \ar[r] & N \times D^j   \, .
}
$$
Then the proposition applied in this case with 
${\cal E} = \Sigma^\infty_N N^+$ (the fiberwise suspension spectrum of
$N^+$) shows that (the stable fiberwise homotopy class of) the map 
$$
f^+ \: P^+ \to N^+
$$
corresponds under the duality isomorphism to the fiberwise collapse
$$
N \times S^j = (N\times D^j)^+ \to P/\!\! /A \, .
$$
This is the fact we are using in the proof of Theorem \ref{uniqueness}.

\end{document}